\documentclass{amsart}

\usepackage{amssymb,amsbsy}

\newtheorem{question}{Question}

\newtheorem{theorem}{Theorem}
\newtheorem{definition}[theorem]{Definition}
\newtheorem*{theorem*}{Theorem}
\newtheorem{proposition}[theorem]{Proposition}
\newtheorem{corollary}[theorem]{Corollary}
\newtheorem*{IVT}{Intermediate Value Theorem}

\newcommand\R{\mathbb R}
\newcommand\vect[1]{\mathbf{#1}}
\newcommand\T{\mathsf{T}}

\newcommand\x{\vect{x}}
\newcommand\y{\vect{y}}
\newcommand\z{\vect{z}}
\newcommand\w{\vect{w}}
\newcommand\0{\vect{0}}
\newcommand\ve{\vect{v}}
\newcommand\be{\vect{b}}
\newcommand\e[1]{\vect{e}_{#1}}
\newcommand\Cone{\mathcal C}

\renewcommand\geq{\geqslant}
\renewcommand\leq{\leqslant}

\begin{document}

\title{The Intermediate Value Theorem \\ for Linear Transformations}
\markright{IVT for Linear Transformations}
\author{Rub\'en A. Mart\'inez-Avenda\~no} 

\maketitle

\begin{abstract}
  If a real-valued function is continuous on a real interval and it takes on two different values, then it will also take any value in between those two, by the Intermediate Value Theorem. It is not immediately clear what would be a natural generalization for functions whose domain and range are in higher-dimensional Euclidean spaces. In this article, we analyze this problem, by first arriving at what we think is the appropriate question to ask, and then restricting to linear transformations. It turns out that the matrices that will satisfy an appropriate version of the Intermediate Value Theorem are the so called {\em monotone} and {\em weakly monotone} matrices, which have applications in numerical approximation of the solutions to systems of linear equations.
\end{abstract}

\section{Introduction.}

In a real analysis course we usually are taught the following very important theorem:

\begin{IVT}
  Let $f: I \subseteq \R \to \R$ be a continuous function on an interval $I$. Assume there are values $y_0, y_1 \in f(I)$, with $y_0 < y_1$. If $ y_0 \leq y \leq y_1$, then $y \in f(I)$. In fact, if $y_0=f(x_0)$ and $y_1=f(x_1)$, there exists $x$ between $x_0$ and $x_1$ such that $f(x)=y$.
  \end{IVT}

Now, what happens if the domain is a subset of $\R^n$?  Suppose that we have a continuous function $f: D \subseteq \R^n \to \R$ and assume $D$ is path-connected (i.e., any two points in $D$ can be joined by a curve contained in $D$). Clearly, if the function $f$ takes the values $c, d \in \R$ , then it must take any value between $c$ and $d$. Indeed, assume $x_0, x_1 \in D$ satisfy $f(x_0)=c$ and $f(x_1)=d$ and, without loss of generality, that $c<d$. Let $\gamma: [0,1] \to D$ be a continuous function (a ``curve'') such that $\gamma(0)=x_0$ and $\gamma(1)=x_1$. Then $f \circ \gamma$ satisfies the hypothesis of the Intermediate Value Theorem and hence, for every $a \in [c,d]$ there exists $t^* \in [0,1]$, such that $(f \circ \gamma)(t^*)=a$. In other words, there exists $x^*:=\gamma(t^*) \in D$ ``between'' $x_0$ and $x_1$ such that $f(x^*)=a$. Thus, for continuous real-valued functions on path-connected subsets of $\R^n$, an Intermediate Value Theorem holds.%
\footnote{In fact, for the more sophisticated reader this is obvious: if $A$ is path-connected, then $A$ is connected. Since $f$ is continuous, then $f(A)$ is a connected set in $\R$ and hence, if $c$ and $d$ are in the interval $f(A)$, so is any point between $c$ and $d$.}

What we would like to explore here, in some sense, is the meaning of ``between'' in the previous paragraph. Yes, $x^*$ is in a path that joins $x_0$ and $x_1$, but it seems an exaggeration to call this behaviour as being ``between'' two points. Fortunately, there is a natural partial order in $\R^n$ which we can use: let $\x=(x_1, x_2, \dots, x_n)$ and $\y=(y_1, y_2, \dots, y_n)$ in $\R^n$; we say that $\x \leq \y$ if $x_j \leq y_j$ for each $j=1, 2, \dots, n$. With this partial order, it is natural to ask the following question:

\begin{question}
Let $f: D \subseteq \R^n \to \R^m$ be a continuous function on a path-connected set $D$. Assume that $\y_0$ and $\y_1$ are in $f(D)$. Given $\y \in \R^m$ with $\y_0 \leq \y \leq \y_1$, is it true that $\y \in f(D)$?
\end{question}

The answer to this question is clearly no. For example, if $f: \R \to \R^2$ is defined as $f(t)=(t,t)$, then clearly $f(0)= (0,0) \leq (2,3) \leq (5,5)=f(5)$, but $(2,3)$ is not in the range of $f$.

What if the domain and codomain coincide? The answer is still no. For example, let $f: \R^2 \to \R^2$ be defined as $f(x,y)=(x^2+y^2, 2 x y)$. Then, since $(0,0) \leq (2,3) \leq (5,4)$, and since $f(0,0)=(0,0)$ and $f(1,2)=(5,4)$, we would expect that $(2,3)$ is also in the image of $f$. But it is not, since the equations $x^2+y^2=2$ and $2 x y = 3$ have no simultaneous solution.

In the latter example, the only vectors $(a,b)$ with $(0,0) \leq (a,b) \leq (5,4)$ that belong to the image of $f$ are those for which $a\geq b$. Thus it seems that characterizing the functions for which the answer to the above question is affirmative might be hard, since the range of these continuous functions might be quite complicated (as opposed to the range of continuous real-valued functions on intervals, which are always intervals). So, perhaps we ought to ask the following question instead.

\begin{question}
  Let $f: D \subseteq \R^n \to \R^m$ be a continuous function on a path-connected set $D$. Let $\x_0$ and $\x_1 \in D$, with $\x_0 \leq \x_1$; set $\y_0:=f(\x_0)$ and $\y_1:=f(\x_1)$ and assume $\y_0 \leq \y_1$. Given $\y \in f(D)$ with $\y_0 \leq \y \leq \y_1$ is there a vector $\x \in D$ with $f(\x)=\y$, and $\x_0 \leq \x \leq \x_1$?
\end{question}

Unfortunately, the answer to this question is also negative. Consider the continuous function $g:\R^2 \to \R^2$ given by $g(x,y)=(x^2-y^2, 2 x y)$. Observe that $g(0,0)=(0,0)$, that $g(7,1)=(48,14)$, that  $(0,0) \leq (5,12) \leq (48,14)$, and that the only solutions of $g(x,y)=(5,12)$, namely $(3,2)$ and $(-3,-2)$, do not satisfy $(0,0) \leq (x,y) \leq (7,1)$. Hence, given $(0,0)$ and $(48,14)$, there exists a point, namely $(5,12)$ {\em in between} $f(0,0)=(0,0)$ and $f(7,1)=(48,14)$ such that none of the points $(x,y)$ satisfying $f(x,y)=(5,12)$ lie in between $(0,0)$ and $(7,1)$. Hence the answer to the question above is also negative. 

Maybe the function above is ``too complicated''?
\footnote{It is not: the astute reader will notice that $g$ is the complex-valued function $g(z)=z^2$.}
What if we try linear functions? In what follows, we intend to answer Question 2 for linear transformations.

\section{Linear Transformations: Questions.}

Let us see what happens if we consider linear functions, so assume $T: \R^n \to \R^m$ is a linear transformation. As usual, we will identify such a linear transformation $T$ with its $m \times n$ matrix $A$ (with respect to the canonical bases) so that $T(\x)= A \x$; from now on we will speak of the matrix $A$ instead of the linear transformation $T$. Elements in Euclidean space will be thought of as column vectors, so that we can multiply them by matrices. When $\vect{v}$ is a row vector, we write $\vect{v}^\T$ to denote its transpose, a column vector.

As we saw in the previous section, Question 1 is false even for linear transformations (see the example right after Question 1). As we mentioned before, one of the problems here is the range of a linear transformation. Let us examine this in more detail. Let $\x=(x_1, x_2, \dots, x_m)^\T$ and $\y=(y_1, y_2, \dots, y_m)^\T$ be vectors in $\R^m$ with $x_k < y_k$ for each $k$. Consider the set
\[
 Z:= \{ \z \in \R^m \mid \x \leq \z \leq \y \}.
\]
In order to answer Question 1, we need to decide if sets like $Z$ are contained in the range of $A$. But, with the conditions imposed on $\x$ and $\y$ above, the set $Z$ contains an open set in $\R^m$. Since the range of a linear tranformation is a vector subspace, and the only subspace of $\R^m$ that contains nontrivial open sets is $\R^m$ itself, it follows that if the set $Z$ is contained in the range of $A$, it must be the case that $A$ is surjective. Hence, unless $A$ is surjective, there will be vectors $\z \in Z$ which are not in the range of $A$. One can get around this issue by requiring some of the coordinates of $\x$ and $\y$ to be equal, but that seems too restrictive.%
\footnote{One can show that for $m \geq 2$, if there exists a unique $k^*$ such that $x_{k^*}=y_{k^*}$, then $Z$ is contained in a nontrivial subspace of $\R^m$ if and only if $x_{k^*}=0$. If $\x$ and $\y$ have $s$ equal coordinates, with $s \geq 2$, then $Z$ is contained in an $(m-s)$--dimensional subspace of $\R^m$ without any further conditions on the values of the coordinates.}
Since we do not want to restrict the choices of $\x \leq \y$, and since if we require $A$ to be surjective then the answer to Question 1 is trivially affirmative, we will forget about Question 1 and instead focus on Question 2, which we rephrase now for linear transformations.

\begin{question}  
  Let $A$ be an $m \times n$ real matrix. For $\x_0$ and $\x_1 \in \R^n$ with $\x_0 \leq \x_1$, set $\y_0:=A \x_0$ and $\y_1:=A\x_1$ and assume that $\y_0 \leq \y_1$. Given $\y$ in the range of $A$ with $\y_0 \leq \y \leq \y_1$, is there a vector $\x \in \R^n$ with $A\x=\y$, and $\x_0 \leq \x \leq \x_1$?
\end{question}

The answer is no. Let
\[
  A=\begin{pmatrix}
    4 & 3 \\
    1 & 1
    \end{pmatrix}.
  \]
  Observe that the matrix $A$ is invertible and hence surjective. Clearly, $A (0,0)^\T = (0,0)^\T$ and $A(3,1)^\T = (15,4)^\T$. Now, consider the vector $(10,3) \in \R^2$: clearly $(0,0) \leq (10,3) \leq (15,4)$ but the only solution of the equation $A \x = (10,3)^\T$ is $\x=(1,2)$ which does not lie between $(0,0)$ and $(3,1)$.

Let us examine another example. Let
\[
  A=\begin{pmatrix}
    1 & 0 & 1 \\
    0 & 1 & 1 
    \end{pmatrix}.
  \]
For $\x_0=(0,0,0)^\T$ and $\x_1=(1,6,3)^\T$, set $\y_0:=A \x_0=(0,0)^\T$ and $\y_1:=A \x_1=(4,9)^\T$. Consider the vector $\y=(3,1)^\T$. Indeed,
  \[
    (0,0) \leq (3,1) \leq (4,9)
  \]
  and since $A$ is surjective there are certainly solutions to the equation $A(x,y,z)^\T=(3,1)^\T$. All solutions to the equation $A(x,y,z)^\T=(3,1)^\T$ are of the form $(3-t, 1-t, t)^\T$, for $t \in \R$. But there is no real value of $t$ for which
\[
  0 \leq 3- t \leq 1, \quad
  0 \leq 1 - t \leq 6, \quad \text{and} \quad
0 \leq t \leq 3.
\]
Hence there is no vector $(x,y,z) \in \R^3$ such that $(0,0,0) \leq (x,y,z) \leq (1,6,3)$ and $(0,0)^\T=A(0,0,0)^\T \leq (3,1)^\T = A (x,y,z)^\T \leq (4,9) = A(1,6,3)^\T$.

This example gives a clue to a better question. Indeed, notice that if we start with the inequality
  \[
    (0,0) \leq (3,1) \leq (4,9),
  \]
  but we consider the vector $(3,8,1)^\T$ (instead of the vector $(1,6,3)^\T$) we also obtain $A(3,8,1)^\T = (4,9)^\T$. In this case we have that  $A(3,1,0)^\T=(3,1)^\T$ and indeed, $(0,0,0) \leq (3,1,0) \leq (3,8,1)$, as desired.

  So, let us rephrase our question once again, by focusing on the vectors in the range of $A$, instead of starting with vectors in the domain.
  
\begin{question}  
  Let $A$ be an $m \times n$ real matrix and let $\y_0$ and $\y_1 \in \R^m$, both in the range of $A$ such that $\y_0 \leq \y_1$. Given $\y$ in the range of $A$ with $\y_0 \leq \y \leq \y_1$, is it true that there exist $\x_0$, $\x$ and $\x_1 \in \R^n$ with $A\x_0=\y_0$, $A\x=\y$, and $A\x_1=\y_1$, and such that $\x_0 \leq \x \leq \x_1$?
\end{question}

In the next section, we attempt to answer Question 4.

\section{Linear Transformations: Answers.}

Before trying to answer Question 4, we need to make a couple of observations. The reader will have noticed that $\x_0=\0$ in all the examples in the previous section. The reason for that should be clear after the following two propositions. In the first one, we show that when considering Question 4 above for only one inequality, we may assume one side of the inequality equals zero. In the second proposition,  we will show that this is in fact equivalent to the inequalities in Question 4.

\begin{proposition}\label{prop:zero_or_not_zero}
  Let $A$ be an $m \times n$ matrix. The following two statements are equivalent.
  \begin{enumerate}
    \item If $A \x_1 \leq A \x_2$, then there exist vectors $\x_1'$ and $\x_2'$ such that $A \x_1= A \x_1'$, $A \x_2 = A \x_2'$, and $\x_1' \leq \x_2'$.
    \item If $\0 \leq A \x$, then there exists a vector $\x'$ such that $A \x= A \x'$,  and $\0 \leq \x'$.
    \end{enumerate}
  \end{proposition}
  \begin{proof}
To prove that (1) implies (2), assume $\0 \leq A \x$. Then, $A \0 \leq A \x$. But (1) implies that there exist vectors $\z_1$ and $\z_2$ such that $A\z_1=A \0=\0$, $A \z_2=A \x$ and $\z_1 \leq \z_2$. Setting $\x':=\z_2-\z_1 \geq \0$ we obtain the desired conclusion, since $A \x'= A(\z_2-\z_1)=A \z_2=A \x$ and $\0 \leq \x'$.

    We now show that (2) implies (1). Assume $A \x_1 \leq A \x_2$. Then we have that $A (\x_2 - \x_1) \geq \0$. But then (2) implies that there exists a vector $\z$ such that $A \z = A (\x_2 - \x_1)$ and $\0 \leq \z$. Then, taking $\x_2':=\x_1+\z$ and  $\x_1':=\x_1$ gives the result, since $\x_2' - \x_1' = \z \geq \0$, $A\x_1'= A \x_1$, and $A\x_2'= A (\x_1 + \z)=  A \x_2$, as desired.
  \end{proof}
    
It is natural to ask, and perhaps not obvious, if the conditions in the above proposition hold not only for two inequalities, but also for three. It turns out that they do hold.
  
  \begin{proposition}
    Let $A$ be an $m \times n$ matrix. The following two statements are equivalent.
    \begin{enumerate}
    \item If $A \x_1 \leq A \x_2$, then there exist vectors $\x_1'$ and $\x_2'$ such that $A \x_1= A \x_1'$, $A \x_2 = A \x_2'$, and $\x_1' \leq \x_2'$.

\item If $A \x_1 \leq A \x_2 \leq A \x_3$, then there exist vectors $\x_1'$, $\x_2'$,  $\x_3'$ such that $A \x_1= A \x_1'$, $A \x_2 = A \x_2'$, $A \x_3= A \x_3'$ and such that $\x_1' \leq \x_2' \leq \x_3'$.
  \end{enumerate}
  \end{proposition}
  \begin{proof}
It is obvious that (2) implies (1). So suppose that (1) holds and assume $A \x_1 \leq A \x_2 \leq A \x_3$. Then we have that
    \[
      A(\x_2-\x_1) \geq \0 \quad \text{ and } \quad A(\x_3-\x_2) \geq \0.
    \]
    Then, Proposition~\ref{prop:zero_or_not_zero} implies that there exists $\z_2\geq \0$ and $\z_3\geq \0$ such that
    \[
      A(\z_2) = A (\x_2 -\x_1) \quad \text{ and } \quad A(\z_3) = A (\x_3 -\x_2).
    \]
Since clearly
    \[
      \x_1 \leq \x_1 + \z_2 \leq \x_1 + \z_2 + \z_3,
      \]
      it follows that if we set $\x_1':=\x_1$, $\x_2':=\x_1 + \z_2$, and $\x_3':=\x_1 + \z_2 + \z_3$, then
\begin{align*}
  A \x_1' &= A \x_1, \\
  A \x_2' & = A( \x_1 + \z_2) = A \x_2, \text{ and } \\
   A \x_3' &= A( \x_1 + \z_2 + \z_3 ) = A \x_2 + A(\x_3-\x_2) = A \x_3,
\end{align*}
which proves the desired result.
    \end{proof}
    
    Hence, to answer question 4, we need to find out for which $m \times n$ matrices $A$ condition (2) of Proposition~\ref{prop:zero_or_not_zero} holds. This leads to the following definition
    \begin{definition}
      Let $A$ be an $m \times n$ matrix. We say that $A$ is {\em weakly monotone} if for all $\x \in \R^n$ such that $\0 \leq A \x$, there exists a vector $\x'$ such that $A \x= A \x'$,  and $\0 \leq \x'$.
      \end{definition}

      It seems that this definition was first given by Adi~Ben-Israel and Thomas~N.E.~Greville in \cite{BIGr}. The full characterization of weakly monotone matrices is not easy and it is beyond the scope of this article (much more information can be found in \cite{BePl,PeSu}). Nevertheless, we will give a partial characterization.

We start with a previous example. Let
      \[
        A=\begin{pmatrix}
          1 & 0 & 1 \\
          0 & 1 & 1 
          \end{pmatrix}.
        \]
To verify that $A$ is weakly monotone, choose $\x=(x_1, x_2, x_3)^\T$ such that $A \x \geq \0$. Then
        \begin{align*}
          x_1 + x_3 &\geq 0\\
          x_2 + x_3 &\geq 0.
         \end{align*}
But clearly we can choose $\x'=(x_1',x_2',x_3')$ with $\x'\geq \0$ with $x_1' + x_3'=x_1 + x_3$ and with $x_2' + x_3'=x_2 + x_3$ (for example, take $x_1':=x_1+ x_3$, $x_2':=x_2+x_3$ and $x_3':=0$) and hence with $A\x = A \x'$. Thus $A$ is indeed weakly monotone.

Before we continue talking about weakly monotone matrices, let us state the following definition which was originally given by Lothar~Collatz \cite{collatz} for square matrices.

\begin{definition}
  Let $A$ be an $m \times n$ real matrix. We say that $A$ is monotone if $\0 \leq A \x $ implies that $\0  \leq \x$.
\end{definition}

Collatz showed that a square matrix is monotone if and only if it is invertible and its inverse is nonnegative.%
\footnote{The anxious reader can try and prove this now: it is not terribly difficult. The patient reader can prove this after seeing the results ahead.}%

Olvi~L.~Mangasarian \cite{mangasarian} extended the above definition to nonsquare matrices and, furthermore, he showed that a matrix is monotone if and only if $A$ has a nonnegative left inverse. Of course, if $A$ has a left inverse, then  $A$ is injective and, in particular, $m \geq n$.

In what follows, we intend to give a partial characterization of weakly monotone matrices, using the techniques in \cite{PeSu, mangasarian}. We  start with a famous lemma, whose proof we include here for the sake of completeness. The proof uses the hyperplane separating theorem for convex subsets of Euclidean space (see, for example, \cite[p.~510]{LuYe}).

\begin{theorem*}[Farkas' Lemma]
  Let $M$ be an $m \times n$ matrix and $\be \in \R^m$. Then one, and only one, of the following statements is true.
  \begin{enumerate}
  \item There exists $\x \in \R^n$ such that $M \x =\be$ and $\x \geq \0$.
  \item There exists $\y \in \R^m$ such that $\y^\T M \geq \0$ and $\y^T \be < 0$.
   \end{enumerate}
    \end{theorem*}
    \begin{proof}
      Assume (1) does not hold. Let $\ve_1, \ve_2, \dots, \ve_n \in \R^m$ be the columns of the matrix $M$. Since (1) does not hold, that means the cone generated by the vectors $\ve_1, \ve_2, \dots, \ve_n$, namely
      \[
        \Cone:=\Bigl\{ \w \in \R^m \mid \w = \sum_{k=1}^m x_k \ve_k, \ x_k \geq 0 \text{ for each } k \Bigl\}
      \]
      does not contain the vector $\be$. Since $\Cone$ is a closed convex set, there exists $\y \in \R^m$ and $\alpha \in \R$ such that the hyperplane $\{ \w \in \R^m \mid \y^\T \w = \alpha \}$ separates $\be$ from $\Cone$. In fact, $\y^\T b < \alpha$ and $\y^T \w > \alpha$ for all $\w \in \Cone$.

      Observe that since $\0 \in \Cone$, we have $\alpha<0$ and hence $\y^\T b < 0$. For each $\lambda >0$ and each $k$, the vector $\lambda \ve_k \in \Cone$. Therefore $\y^\T \lambda \ve_k > \alpha$; equivalently, $\y^\T \ve_k > \frac{\alpha}{\lambda}$. Since $\lambda>0$ can be arbitrarily large (and $\alpha < 0$), this means that $\y^\T \ve_k \geq 0$ for all $k$. But this is equivalent to $\y^T M \geq \0$.  That is, (2) holds.

      Now, let us see that (1) and (2) cannot both be true. By contradiction, assume there exists $\x \in \R^n$ such that $M \x =\be$ and $\x \geq \0$ and there exists $\y \in \R^m$ such that $\y^\T M \geq \0$ and $\y^T \be < 0$. But then on one hand
      \[
        \y^\T M \x = \y^T \be < 0,
      \]
      and on the other $\y^\T M \x \geq 0$, since $\y^\T M \geq \0 $ and  $\x \geq \0$. This is a contradiction, and the result follows.      
      \end{proof}

This allows us to give a sufficient condition for a matrix to be weakly monotone. Recall an $m \times n$ matrix $M$ is called {\em nonnegative}, denoted by $M \geq 0$, if all of its entries are nonnegative. Here $I_n$ and $I_m$ denote the identity matrices of size $n$ and $m$, respectively.
      
      \begin{proposition}\label{prop:suff}
        Let $A$ be an $m \times n$ matrix.
        \begin{itemize}
        \item If $m \geq n$, assume there exists an $n \times m$ matrix $B$ with $B A = I_n$ such that $B \geq 0$.
        \item If $m \leq n$, assume there exists an $n \times m$ matrix $B$ with $A B = I_m$ such that $B \geq 0$.
        \end{itemize}
        Then for all $\x \in \R^n$ we have that $A \x \geq \0$ implies that there exists $\x' \in \R^m$ such that $\x' \geq \0$ and $A \x = A \x'$.
        \end{proposition}
\begin{proof}
Assume first that $m \geq n$. Let $B$ be an $n \times m$ matrix such that $B A=I_n$ and $B \geq 0$. Let $\x \in \R^n$ such that $A \x \geq \0$, then, since $B \geq 0$, we have $B A \x \geq \0$. But then $\x \geq \0$, so taking $\x'=\x$ finishes this case. 

Assume now that $m \leq n$. Let $B$ be an $n \times m$ matrix such that $A B =I_m$ and $B \geq 0$. Let $\x \in \R^n$ and assume $\be:=A \x \geq \0$. If there is no $\x' \geq \0$ such that $A \x'=\be$, then Farkas' Lemma implies that there exists $\y \in \R^m$ such that $\y^T A \geq \0$ and $\y^\T \be < 0$. Since $B \geq 0$ this implies that $\y^T A B \geq \0$; hence $\y^\T \geq \0$. Since $\be \geq \0$, this implies  $\y^\T \be \geq 0$, which is a contradiction. This finishes the proof.
\end{proof}      

Unfortunately, as we will see later, the above conditions are not necessary. Nevertheless, they are if we give a further condition on $A$.

\begin{proposition}\label{prop:nec}
  Let $A$ be an $m \times n$ matrix of full rank. Assume that for all $\x \in \R^n$ with $A \x \geq \0$ there exists $\x' \in \R^m$ such that $\x' \geq \0$ and $A \x = A \x'$. Then,
  \begin{itemize}
  \item If $m \geq n$, there exists an $n \times m$ matrix $B$ with $B A = I_n$ such that $B \geq 0$.
  \item If $m \leq n$, there exists an $n \times m$ matrix $B$ with $A B = I_m$ such that $B \geq 0$.
    \end{itemize}
  \end{proposition}
  \begin{proof}
    Assume first that $m \geq n$. Since $A$ is of full rank, we have that the rank of $A$ is $n$. This means that $A$ is injective.
    Let $\e{k} \in \R^n$ be the $k$-th canonical vector $\e{k}=(0, 0, \dots, 0, 1, 0, \dots, 0)^\T$, with the number $1$ in the $k$-th coordinate. We claim that, for each $k=1, 2, \dots n$, the equation $A^\T \x = \e{k}$ has a solution $\x^{k} \geq \0$, $\x^{k}\in \R^m$; for if not, Farkas' Lemma would imply that there exists $\y \in \R^n$ such that $\y^\T A^\T \geq 0$ and $\y^\T \e{k} <0$. But this is equivalent to saying that $A \y \geq 0$ and the $k$-th entry of $\y$ is negative. But since $A$ is injective, if we set $\be:=A \y$, the only solution to the equation $A \x = \be$ is $\x=\y$. Since the $k$-th entry of $\y$ is negative, this contradicts the hypothesis. Therefore, for each $k$ we have that there exists $\x^{k} \geq \0$ such that
    \[
      (\x^{k})^\T A = \e{k}^T.
    \]
    Taking $B$ to be the $n \times m$ matrix whose $k$-th row is the vector $(\x^{k})^\T$, gives that $B A = I_n$. Since $\x^{k} \geq \0$, it follows that $B \geq 0$. 

Now assume $m \leq n$. Since $A$ is of full rank, we have that rank of $A$ is $m$. Hence, for each $k=1, 2, \dots m$, the equation $A \x = \e{k}$ has a solution $\x \in \R^n$. In fact, by the hypothesis, one can choose $\x^{k} \geq 0$, such that $A \x^{k}=\e{k}$. Taking $B$ to be the $n \times m$ matrix whose $k$-th column is the vector $\x^{k}$, gives that $A B = I_m$. Since $\x^{k} \geq \0$, it follows that $B \geq 0$, which finishes the proof.
  \end{proof}

There is a case in which Propositions~\ref{prop:suff} and \ref{prop:nec} can be combined into an easy necessary and sufficient condition as the next corollary states.
  
  \begin{corollary}
    Let $A$ be an $1 \times n$ nonzero matrix. The following conditions are equivalent:
    \begin{enumerate}
    \item For all $\x \in \R^n$ if we have that $A \x \geq 0$ then there exists $\x' \in \R^m$ such that $\x' \geq 0$ and $A \x = A \x'$.
    \item There exists an $n \times 1$ matrix $B$ with $A B = 1$ such that $B \geq 0$.
    \item The matrix $A$ has at least one positive entry.
    \end{enumerate}
  \end{corollary}
  \begin{proof}
    By Propositions \ref{prop:suff} and \ref {prop:nec}, it is clear that (1) is equivalent to (2), since $A$ is of full rank. That (3) implies (2) is clear: if the $k$-th entry of $A$, say $a_k$ is positive, take $B$ to be the $n \times 1$ matrix with $\frac{1}{a_k}$ on the $k$-th entry and zeros elsewhere. Now, assume (2) holds but (3) is not true; that is, all entries of $A$ are nonpositive. But then, multiplying $A$ by any $n \times 1$ matrix $B$ with $B\geq 0$ must give a nonpositive number, contradicting (2).
    \end{proof}

    Observe that the corollary above implies that the answer to Question 4 is affirmative, for $m=1$, if and only if $A$ has a positive entry. This has a nice geometric interpretation: if we think of the matrix $A$ as a vector in $\R^n$, in order to have an $n \times 1$ matrix $B$ with $A B = 1$ and with $B \geq 0$, it is necessary and sufficient that the dot product of the ``vector'' $A$ with a vector in the first orthant%
\footnote{The first orthant of $\R^n$ is the set of vectors $\x \in \R^n$ with $\x \geq \0$.} 
of $\R^n$ be positive; i.e., its angle with a vector in the first orthant of $\R^n$ is less than $90$ degrees. This is possible as long as the ``vector'' $A$ is not in the opposite orthant of the first orthant; i.e., at least one entry of $A$ is positive (compare with \cite[Theorem 1]{mangasarian}).

Now there are examples of weakly monotone matrices not of full rank, and hence we cannot apply Proposition~\ref{prop:nec} to them. For example, if we let
  \[
    A=\begin{pmatrix}
      1 & 0 & 1 \\
      0 & 1 & 1 \\
      0 & 0 & 0
      \end{pmatrix},
    \]
then $A$ is of rank $2$. An argument similar to the one before Definition 2 shows that $A$ is weakly monotone, but $A$ cannot have a left or right inverse (in particular, $A$ is not monotone).

Now, to finish this section, let us briefly explore what can we do to answer Question 4 if the matrix $A$ is not of full rank. So assume that we are given an $m \times n$ matrix $A$ and a vector $\be \geq \0$ such that the equation $A \x = \be$ has a solution and we want to find $\x' \geq 0$ such that $A \x'=\be$. The natural course of action is to find all solutions to $A \x = \be$ and check if we can find a nonnegative one. The usual procedure is to find an invertible $m \times m$ matrix $Q$ such that
\[
Q A P = \left( \begin{array}{c|c}
               I_k & S \\
               \hline
               0 & 0 
\end{array}\right),
\]
where $P$ is an $n \times n$ permutation matrix, $I_k$ is the $k \times k$ identity matrix and $S$ is a $k \times (n-k)$ matrix, and the zero matrices, denoted by $0$, are of appropriate sizes (here $k$ is the rank of $A$, of course). We should point out that in some cases, there may be several choices of matrices $Q$ and $P$ (and $S$).

If $Q$ is a nonnegative matrix, then all works out and hence $A$ is weakly monotone. Indeed, left multiplying by $Q$ both sides of the equation $A \x = \be$ we obtain the equivalent equation
\[
  Q A \x = Q \be, \quad \text{ or, equivalently } \quad   Q A P P^\T \x = Q \be,
\]
where we are using that $P P^T=I_n$, since $P$ is a permutation matrix. Since we know the equation has a solution, it must be that $Q \be = ( \be' \ | \ \0 )^\T$, for some $\be' \in \R^k$ and $\0 \in \R^{m-k}$. Furthermore, $\be' \geq 0$, since $Q$ and $\be$ are nonnegative. But this implies that if we set $P^\T \x = (\y \ | \  \z)^T$, where $\y \in \R^k$ and $\z \in \R^{n-k}$ we obtain
\[
  (\y + S \z \ | \  \0 )^\T = Q \be = ( \be' \ | \ \0 )^\T,
\]
and hence $\y = \be'-S\z$. Setting $\z=\0$, gives that $P^\T \x = (\be' \ | \  \0)^T$, and hence $\x' = P (\be' \ | \  0)^T$  is the sought-after nonnegative solution.

Let is illustrate this procedure with an example. Let
\[
  A=\begin{pmatrix}
    1 & -1 &  4  \\
   -1 &  2 & -3 \\
    1 & -2 &  3
  \end{pmatrix},
\]
this leads to $A (-3,0,1)^\T=(1,0,0)^\T=:\be$.
Then, we can find, by elementary row operations, the matrix
\[
  Q=\begin{pmatrix}
    2 & 1 & 0 \\
    1 & 1 & 0 \\
    0 & 1 & 1 
    \end{pmatrix}
  \]
  such that $QA$ is in reduced row echelon form. But then the equation $A \x = \be$ is equivalent to $Q A \x = Q \be$. This is the equation
  \[
 \begin{pmatrix}
      1 & 0 & 5 \\
      0 & 1 & 1 \\
      0 & 0 & 0      
      \end{pmatrix} \begin{pmatrix} y_1 \\ y_2 \\ z_1 \end{pmatrix} = \begin{pmatrix} 2 \\ 1 \\ 0  \end{pmatrix},
    \]
and hence taking $z_1=0$, we obtain a nonnegative solution $\x'=(2,1,0)^\T$ to the equation $A \x = \be$.

\section{Final remarks.}

To summarize what we have done in this article, if $A$ is a weakly monotone matrix, the answer to Question 4 is affirmative. As we mentioned before, characterizing which matrices are weakly monotone is not a simple problem (again, see \cite{PeSu}). Nevertheless, checking if a single matrix $A$ is weakly monotone or not is a straightforward (but perhaps not computationally easy) problem, as we mentioned in the previous section.

Surprinsingly perhaps, the answer to Question 3 is easier. The $m \times n$ matrix $A$ is monotone if and only if it has a left inverse (and hence $m \geq n$), and if this occurs then the answer to Question 3 is affirmative. Fewer matrices will have this stronger property.

We should mention that the answers to Questions 3 and 4 above have important applications. Collatz~\cite{collatz} introduced the concept of mononote matrices in his study of the approximation of solutions of systems of linear equations. Indeed, for a nonsingular square matrix $A$, if we need to solve the equation $A \x = \be$, and we have found vectors $\ve_1 \leq \ve_2$ such that $A \ve_1 \leq \be  \leq A \ve_2$, it seems natural to expect that the solution $\x$ might be found ``between'' $\ve_1$ and $\ve_2$; i.e., $\ve_1 \leq \x \leq \ve_2$. Collatz characterized those matrices for which this occurs: monotone matrices.

It is not surprising then than the concept of monotonicity is crucial in the study of iterative methods for solving systems of linear equations. Richard~S.~Varga~\cite{varga} introduced the concept of regular splitting of a matrix $A$, which includes the monotonicity of the matrix $A$, to find necessary conditions for the convergence of certain iterative methods for approximating the solutions to $A \x = \be$.  See \cite{BIGr,BePl,MePl} and the references therein for more information on these topics, including some applications to linear programming.

What is it about linear maps that allowed us to partially solve Question 4? Observe that in Proposition~\ref{prop:suff}, in the case $m \geq n$ we used the existence of a left inverse that sends nonnegative vectors to nonnegative vectors. In general, for a function $f: \R^m \to \R^n$, if there exists a function $g:\R^n \to \R^m$ which sends nonnegative vectors to nonnegative vectors, and such that $g \circ f$ is the identity on $\R^m$, then we would have that if $\0 \leq f( \x )$, then $\0 \leq g(f(\x)) = \x$, so condition (2) of Proposition~\ref{prop:zero_or_not_zero} is satisfied (of course, since $f$ is not necessarily linear, we do not have the equivalence in Proposition~\ref{prop:zero_or_not_zero}). When is there a left inverse $g$ to a function $f$, such that $g$ maps the first orthant into itself? Perhaps nonlinear versions of Farkas' Lemma (e.g, \cite{Craven, Zalinescu, Glover}) might give a partial answer. We leave this question open for the interested reader.

Lastly, we should point out that there exist other generalizations of the Intermediate Value Theorem in several dimensions; e.g., the Poincar\'e-Miranda theorem (see, for example, \cite{kulpa}).

\vskip0.5cm

{\bf Acknowledgments.}
  {The author wishes to thank Levent Ulku, who approached the author with a question that led to the topic of this paper. Also, the author thanks Blair Madore, Edgar Possani and the anonymous reviewers for their comments on the first version of this manuscript, which significantly improved the exposition.
\\
This work is partially supported by the Asociaci\'on Mexicana de Cultura A.C.
\\
2020 Mathematics Subject Classification: Primary 15A24, Secondary 15A39, 15B48.

\end{document}